\documentclass[english,reqno]{amsart}

\usepackage{amsmath, appendix, ulem}
\usepackage[nobysame]{amsrefs}
\usepackage{amssymb, color}
\usepackage[margin=1in]{geometry}

\usepackage{mathrsfs}
\usepackage{graphicx}
\usepackage{subfig}
\usepackage{float}
\usepackage{epsf}
\usepackage{hyperref}
\usepackage{titletoc}

\usepackage{CJKutf8}

\usepackage{MnSymbol}

\graphicspath{{../Figures/}}

\numberwithin{equation}{section}

\newtheorem{thm}{Theorem}[section]
\newtheorem{prop}[thm]{Proposition}

\theoremstyle{remark}

\renewcommand{\tilde}{\widetilde}

\renewcommand{\bar}{\overline}

\newcommand{\eq}{\mathrm{eq}}
\newcommand{\nn}{\nonumber}
\newcommand{\R}{{\mathbb R}}

\newcommand{\del}{\partial}

\newcommand{\ds}{\, {\rm d} s}
\newcommand{\dt}{ \, {\rm d} t}

\newcommand{\dx}{ \, {\rm d} x}
\newcommand{\dy}{ \, {\rm d} y}
\newcommand{\dz}{ \, {\rm d} z}

\newcommand{\dmu}{\, {\rm d} \mu}

\newcommand{\One}{\boldsymbol{1}}

\newcommand{\Eps}{\epsilon}
\newcommand{\Ni}{\noindent}

\newcommand{\CalM}{{\mathcal{M}}}

\newcommand{\domega}{\, {\rm d}\omega}

\newcommand{\abs}[1]{\left\lvert#1\right\rvert}
\newcommand{\norm}[1]{\left\lVert#1 \, \right\rVert}

\newcommand{\dvint}[1]{\left\llangle#1\right\rrangle}

\newcommand{\vpran}[1]{\left(#1\right)}

\begin{document}

\title{Unique reconstruction of the heat-reflection indices at solid interfaces}

\author{Qin Li} 
\address{Mathematics Department, University of Wisconsin-Madison, 480 Lincoln Dr., Madison, WI 53705 USA.}
\email{qinli@math.wisc.edu}
\author{Weiran Sun}
\address{Department of Mathematics, Simon Fraser University, 8888 University Dr., Burnaby, BC V5A 1S6, Canada}
\email{weirans@sfu.ca}

\date{\today}

\thanks{The research of Q.L. is supported in part by NSF under DMS-1750488, and ONR-N000142112140. The research of W.S. is supported in part by NSERC Discovery Grant R611624.}

\begin{abstract}
We show the unique reconstruction of the heat-reflection coefficients in a phonon transport equation. This is a mathematical model used to characterize the dynamics of heat-conducting phonons in multiple layers of media, commonly composed of metals and silicon. In experiments, the heat-reflection indices are inferred by measuring the temperature at the surface of the exterior metal after applying heat sources. In this article, we rigorously justify the unique reconstruction of these indices by using such procedures. 
\end{abstract} 
%
\maketitle

\section{Introduction}
How does heat propagate through the interface between two solids? Traditionally the heat equation was used to describe the process. More recently, it is discovered, from the first principle, that phonon transport equations containing microscopic information give a more accurate description~\cite{10.1115/1.4002028}. The heat conduction through interfaces is determined by the reflection/transmission coefficients in the phonon transport equation~\cite{PhysRevB_Hadjiconstantinou}. These coefficients vary among different materials. Non-intrusive experiments taken on the surface of the materials have been used to reconstruct them~\cite{PhysRevB_Minnich,PhysRevX.8.041004}. A natural question is to decide how to design the experiment and which data to measure so that a unique reconstruction of the coefficient is guaranteed.

In this paper, we study the experimental setup used in~\cite{PhysRevB.95.205423}, where two solids are placed side by side with one being a metal (aluminum e.g.) and the other being silicon. Heat is injected at the surface of the metal and propagates to the silicon, resulting in a temperature increase in both materials. One then measures the temperature on the surface of the metal (as a function of time) to deduce the reflection indices at the interface of the two solids. 
To numerically reconstruct the parameters, the authors in~\cite{PhysRevB.95.205423} relied on a minimization formulation. In this paper, we rigorously prove that by carefully designing the injecting phonons, this experimental setup and its measurement will uniquely recover the reflection indices. Hence this experimental method is supported by a solid theoretical foundation. The main idea of our proof is built upon the classical singular decomposition approach. In particular, we will choose the incoming data to be singular in time, speed and frequency. The observation is that along the trajectory of the heat transport, such singularity is almost preserved. As a consequence, the majority of the information bounced back at the interface and captured at the surface in the measurement contains specific values of the reflection coefficients.




Phonon transport equations belong to the larger framework of kinetic theory, which characterizes the dynamics of particles from a statistical perspective~\cite{Cercignani}. It is common that parameters in kinetic equations that encode the intrinsic information of the media or particles are unknown. Thus in many real-world applications, kinetic theory is presented in its inverse facet: measurements need to be taken to deduce these parameters. The problem studied in this article is one such example. Many other inverse kinetic problems have been studied. Examples are the Vlasov-Poisson system for electrons in semiconductors~\cite{Gamba} and the radiative transfer equation for photon dynamics~\cite{Bal_Chung_Schotland,Ammari_Jin,Ren_fPAT,Arridge_Schotland09}. In~\cite{LiSun_2020}, the authors studied the general use of kinetic tools (singular decomposition with averaging lemma) in inverse kinetic theory and in~\cite{GLN20}, some numerical studies have been performed.

Most of the approaches taken to derive parameters in inverse problems trace back to the classical singular decomposition method. The main idea of this method is to single out the component of the solution that best reflects the configuration of the unknown parameter. It is particularly useful for kinetic equations since the loss term in the collision operator, when combined with the transport term, typically resembles the attenuated X-ray transform. Meanwhile, the nonlocal terms, or the gain terms in many models, can be made negligible by suitably designing the initial/boundary data. Such procedure reduces the complicated integro-differential kinetic equations to the well-understood X-ray problem. The earlier exploration of this approach was  in~\cite{CS98,CS2}, and it has been extended to kinetic problems in other scenarios~\cite{SU2d,Stefanov_Tamasan,BalMonard_time_harmonic,Bal14,Wang1999,LLU18}.

The rest of the paper is laid out as follows.  In Section~\ref{sec:setup} we present the equation in its general form. The physical setup and the problem's mathematical description are also discussed in this section. In Section~\ref{sec:APriori} we show some a-priori $L^p$ estimates to the solutions. The reconstruction will be detailed in Section~\ref{sec:theory}. 

\section{Mathematical setup}\label{sec:setup}
In this section we explain the model and experimental setup used in~\cite{PhysRevB.95.205423}. First, we introduce the phonon transport equation and its linearization around the room temperature. Let $f(t,x,\mu,\omega)$ be the energy distribution function of phonon particles at time $t>0$, location $x \in \R$, frequency $\omega \in \R^+$,  with velocity $\mu \in [-1, 1]$. The problem has been reduced to one dimension since the experiments are conducted between large parallel planes, thus plane symmetry is applicable. In this setting, the phonon transport equation has the form
\begin{align} \label{eq:main-nonlinear}  
   \del_t f + \mu v_\omega \del_x f = \frac{1}{\tau(\omega)} \vpran{M_f - f},
\end{align}
where $\tau(\omega)$ is the relaxation time and $v_\omega$ is the group velocity. The term on the left-hand side of~\eqref{eq:main-nonlinear} describes the standard phonon transport and the term on the right characterizes the interaction of phonons with the underlying heat bath. The right-hand side is usually called the BGK term and $M_f$ is an equilibrium state. In~\eqref{eq:main-nonlinear}, $M_f$ is related to the Bose-Einstein distribution and is defined as
\begin{align} \label{def:M-f}
   M_f = \frac{\omega}{e^{\frac{\hbar \omega}{k_0 T[f]}} - 1}, 
\end{align}
where $\hbar$ is the Planck constant and $k_0$ the Boltzmann constant. The temperature $T[f]$ is a functional of $f$. It is defined through the energy conservation as follows: integrate the right-hand side and set it to $0$ to get
\begin{align} \label{eq:conserv}
  \int_0^\infty \int_{-1}^1 \frac{1}{\tau(\omega)} \vpran{M_f - f} \dmu \domega = 0.
\end{align}
Inserting the definition of $M_f$ into~\eqref{eq:conserv} gives
\begin{align} \label{def:h}
   h(T) := \int_0^\infty \frac{1}{\tau(\omega)}    \frac{\omega}{e^{\frac{\hbar \omega}{k_0 T}} - 1} \domega = \int_{-1}^1 \int_0^\infty \frac{1}{\tau(\omega)} f(t, x, \mu, \omega) \domega \dmu\,,
\end{align}
which implicitly defines $T$ as a functional of $f$. Note that $h$ is invertible since it is strictly increasing in $T$.

Since experiments are usually conducted at the room temperature $T_\eq$, the linearization of $T[f]$ around $T_\eq$ is a valid approximation~\cite{PhysRevB.95.205423}. To derive the linearized equation, let $\tilde{f} = f-M_\eq$ where $M_\eq$ is defined in~\eqref{def:M-f} with $T=T_\eq$. Then the BGK term becomes
\begin{align*}
M_f-f & =  M_{M_\eq + \tilde f} - (M_{M_{\eq}} - \tilde f) \approx M'|_{T_{\eq}} \vpran{T[F_{\eq} + \tilde f] - T_\eq} - \tilde f,
\end{align*}
where $M'$ stands for $\partial_TM$. To compute $T[F_{\eq} + \tilde f] - T_\eq$, we use the definition of $h$ in~\eqref{def:h} and obtain
\begin{equation*}
M_f-f\approx \frac{M'|_{T_{\eq}}}{h'|_{T_{\eq}}} \int_{-1}^1 \int_0^\infty \frac{1}{\tau(\omega)} \tilde f(t, x, \mu, \omega) \domega \dmu - \tilde f\,.
\end{equation*}
The derivatives of $M$ and $h$ can be computed explicitly and we get
\begin{align*}
   \frac{M'|_{T_{\eq}}}{h'|_{T_{\eq}}}= c_0 M_{\eq}^2 e^{\frac{\hbar \omega}{k_0 T_{\eq}}}, \qquad\text{with}\qquad   c_0 = \frac{\hbar}{k_0 T_{eq}^2 \, h'|_{T_{eq}}}\,.
\end{align*}
Denote 
\begin{align}\label{eqn:bracket}
   \xi(\omega) = \frac{M'|_{T_{eq}}}{h'|_{T_{eq}}} \geq 0\,,
\qquad \qquad
 \dvint{\tilde f} = \int_{-1}^1 \int_0^\infty \frac{1}{\tau(\omega)} \tilde f(t, x, \mu, \omega) \domega \dmu\,.
\end{align}
Then $\xi \in L^\infty(0, \infty)$ and by the conservation law, it satisfies the normalization
\begin{align*}
   \int_0^\infty \frac{\xi(\omega)}{\tau(\omega)} \domega = 1.
\end{align*}
Moreover, by the exponential decay of $M_{\eq}$, we also have
\begin{align} \label{bound:xi-1}
   \int_0^\infty \xi^q (\omega) \domega < C_q < \infty\,,
\qquad\forall q>0\,.
\end{align}
Summarizing the calculation above, we obtain the linearized operator as
\begin{align*}
   L \tilde f = \frac{1}{\tau(\omega)} \vpran{\xi(\omega) \dvint{\tilde f} -  \tilde f}.
\end{align*}
Drop the tilde signs for the simplicity of the notation. Then the linearized phonon transport equation is
\begin{align} \label{eq:linear-basic}
   \del_t f + \mu v_\omega \del_x f = L f 
= \frac{1}{\tau(\omega)} \vpran{\xi(\omega) \dvint{f} -  f}.
\end{align}

The particular experiment places two layers of solid materials side by side. Each material is modelled by a linear phonon transport equation in the interior, and the two equations are coupled at the interface with reflection and transmission between these two solids. The coupled system writes
\begin{align} \label{sys:full-f}
\begin{split}
 &  \del_t f + \mu v_\omega \del_x f 
= \frac{1}{\tau(\omega)} \vpran{\xi(\omega) \dvint{f} -  f}\,,\hspace{1.35cm} x\in[0,1]\,,
\\
& f |_{x=0} = \phi(t, \mu, \omega), \hspace{4.4cm} \mu > 0\,,
\\
&  f|_{x=1} = \eta_1(\omega) f(1, -\mu, \omega) + \zeta_1(\omega) g(1, \mu, \omega), \qquad \mu < 0\,,
\\
& f|_{t=0} = 0
\end{split}
\end{align}
and
\begin{align} \label{sys:full-g}
\begin{split}
 &  \del_t g + \mu v_\omega \del_x g 
= \frac{1}{\tau(\omega)} \vpran{\xi(\omega) \dvint{g} -  g} \,,\hspace{1.35cm} x\in[1,L]\,,
\\
& g |_{x=1} = \eta_2(\omega) f(1, \mu, \omega) + \zeta_2(\omega) g(1, -\mu, \omega), \hspace{1.5cm} \mu > 0\,,
\\
&  g |_{x=L} = \alpha_0 \xi(\omega) \int_0^1 \int_0^\infty \mu \, v_\omega \, g(t, L, \mu, \omega) \domega \dmu, \qquad \mu < 0\,,
\\
& g|_{t=0} = 0\,.
\end{split}
\end{align}
In this model, $f$ and $g$ are distribution functions for phonons in the two materials respectively. The metal occupies the spatial interval $[0,1]$ and the silicon is in $[1,L]$. Both distribution functions have zero initial data. The boundary condition for $f$ is of incoming-type at $x=0$ and reflective at the interface $x=1$, and the boundary condition for $g$ at $x=L$ is diffusive -- the phonons are bounced back with the equilibrium profile. The coefficient $\alpha_0$ is the normalization constant which ensures the energy flux at $x=L$ is zero. This gives
\begin{align*}
   \alpha_0 = \frac{1}{\int_0^1 \int_0^\infty \mu \, v_\omega \, \xi(\omega) \domega \dmu}.
\end{align*}
The most interesting physics takes place at $x=1$ where the two materials meet. Heat propagates through the solid interface according to the transmission and reflection coefficients $\eta_1, \eta_2, \zeta_1, \zeta_2$. These four coefficients satisfy the relation
\begin{align} \label{assump:coeff}
   \eta_1 + \eta_2 = 1, 
\qquad
  \zeta_1 + \zeta_2 = 1,
\qquad
  \eta_1 + \gamma_0 \zeta_1 = 1
\end{align}
for some constant $\gamma_0 > 0$.
These conditions guarantee the conservation of energy across the interface and that $(\xi, \gamma_0 \xi)$ is an equilibrium solution to the coupled system~\eqref{sys:full-f}-\eqref{sys:full-g} (except for the initial data). Because of~\eqref{assump:coeff}, one only needs to reconstruct the reflection coefficient $\eta_1$ and all others follow by simple algebra.

A (very) brief summary of our main result is
\begin{thm} [Main Theorem]
By proper choices of the incoming data and the measurement function, the heat-reflection coefficients $(\eta_1, \eta_2, \zeta_1, \zeta_2)$ can be uniquely and explicitly reconstructed. 
\end{thm}

\section{A priori estimates}\label{sec:APriori}

In this section we prove some a priori $L^p$ estimates for the solution $(f, g)$. First we show the bounds for $p = 1$ and $p = \infty$. The general bound for $p \in (1, \infty)$ follows from the Riesz-Thorin interpolation theorem.

In what follows we will use the notations $L^p_{\mu, \omega}$ and $L^p(\!\dmu \domega)$ interchangeably to denote 
the usual $L^p$-space in the $(\mu, \omega)$ variables for $1 \leq p \leq \infty$. If there is a weight $h$ in the measure then we denote the space as $L^p(h \dmu \domega)$. Similar rules apply to $L^p_{x, \mu, \omega}$ and its weighted version. 

\begin{prop} \label{prop:f-L-p}
Suppose $\phi$ is the incoming data for the system~\eqref{sys:full-f}-\eqref{sys:full-g} and $(f, g)$ is the solution. 
\smallskip

\Ni (a) If $\phi \in L^1_{\mu, \omega}$, then for each $t \geq 0$, 
\begin{align*}
   \norm{f(t, \cdot, \cdot, \cdot)}_{L^1_{x, \mu, \omega}}
+ \norm{g(t, \cdot, \cdot, \cdot)}_{L^1_{x, \mu, \omega}}
\leq
   \norm{\phi}_{L^1(0, t; L^1_{\mu, \omega})}. 
\end{align*}

\Ni (b) Suppose there exists a constant $m_0 > 0$ such that
\begin{align*}
   0 \leq \phi \leq m_0 \xi(\omega). 
\end{align*}
Then for each $t \geq 0$, 
\begin{align*}
  0 \leq f(t, x, v, \omega) \leq m_0 \xi(\omega), 
\qquad
  0 \leq g(t, x, v, \omega) \leq \gamma_0 m_0 \xi(\omega).
\end{align*}

\Ni (c) For any $p \in (1, \infty)$, if $\phi \in L^p(0, t; L^p(\xi^{1-p} \dmu\domega))$,  then
\begin{align} \label{bound:L-p}
   \norm{(f, g)(t, \cdot, \cdot, \cdot)}_{L^p(\xi^{1-p} \dx \dmu \domega)}
\leq
   (1 + \gamma_0)^{1/p} \norm{\phi}_{L^p(0, t; L^p(\xi^{1-p} \dmu\domega))}. 
\end{align}

\end{prop}

Most of the proof follows from classical kinetic theory techniques which can be traced back to~\cite{DL}. We show the details below for the completeness of the work. 
\begin{proof}
(a) The absolute values $\abs{f}$ and $\abs{g}$ satisfy the system
\begin{align} \label{sys:full-f-abs}
\begin{split}
 &  \del_t \abs{f} + \mu v_\omega \del_x \abs{f}
\leq \frac{1}{\tau(\omega)} \vpran{\xi(\omega) \dvint{\abs{f}} -  \abs{f}}, 
\\
& \abs{f} |_{x=0} \leq \abs{\phi(t, \mu, \omega)}, \hspace{3.6cm} \mu > 0,  
\\
&  \abs{f} |_{x=1} \leq \eta_1 \abs{f(1, -\mu, \omega)} + \zeta_1 \abs{g(1, \mu, \omega)}, \qquad \mu < 0, 
\\
& \abs{f}|_{t=0} = 0
\end{split}
\end{align}
and
\begin{align} \label{sys:full-g-abs}
\begin{split}
 &  \del_t \abs{g} + \mu v_\omega \del_x \abs{g} 
\leq \frac{1}{\tau(\omega)} \vpran{\xi(\omega) \dvint{\abs{g}} -  \abs{g}}, 
\\
& \abs{g} |_{x=1} \leq \eta_2 \abs{f(1, \mu, \omega)} + \zeta_2 \abs{g(1, -\mu, \omega)}, \hspace{2.9cm} \mu > 0,   
\\
&  \abs{g} |_{x=L} \leq \alpha_0 \xi(\omega) \int_0^1 \int_0^\infty \mu' \, v_{\omega'} \, \abs{g(t, 1, \mu', \omega')} \domega' \dmu', \qquad \mu < 0,  
\\
& \abs{g} |_{t=0} = 0.  
\end{split}
\end{align}
Noting that if we integrate $\mu v_\omega |f|$ and $\mu v_\omega |g|$ in $(x, \mu, \omega)$ and use the boundary conditions at the interface:
\begin{align*}
&   \int_0^\infty \int_{-1}^1 \mu v_\omega |f| \big|_{x=1} \dmu \domega 
   - \int_0^\infty \int_{-1}^1 \mu v_\omega |f| \big|_{x=0} \dmu \domega
\\
& \hspace{0.3cm} \geq
   \int_0^\infty \int_0^1 (1 - \eta_1) \mu v_\omega |f|(1, \mu, \omega) \dmu \domega
  + \int_0^\infty\int_{-1}^0 \zeta_1 \mu v_\omega  |g|(1, \mu, \omega) \dmu \domega
  - \int_0^\infty \int_0^1  \mu v_\omega |\phi_0| \dmu \domega
\end{align*}
and
\begin{align*}
&  \int_0^\infty \int_{-1}^1 \mu v_\omega |g| \big|_{x=L} \dmu \domega
   - \int_0^\infty \int_{-1}^1 \mu v_\omega |g| \big|_{x=1} \dmu \domega
\\
& \hspace{1cm} \geq
  -\int_0^\infty \int_{-1}^0 (1 - \zeta_2) \mu v_\omega |g|(1, \mu, \omega) \dmu \domega
  - \int_0^\infty \int_0^1 \eta_2 \mu v_\omega  |f|(1, \mu, \omega) \dmu \domega. 
\end{align*}
Adding the two sets of equations and utilizing the above-calculated boundedness gives
\begin{align*}
   \frac{\rm d}{\dt} 
   \vpran{\norm{f}_{L^1_{x, \mu, \omega}} + \norm{g}_{L^1_{x, \mu, \omega}}}
\leq
  \int_0^1 \int_0^\infty \mu v_\omega |\phi_0(t, \mu, \omega)| \dmu \domega.
\end{align*}
Hence, for any $t \geq 0$, it holds that
\begin{align*}
  \norm{f}_{L^\infty(0, t; L^1_{x, \mu, \omega})} 
  + \norm{g}_{L^\infty(0, t; L^1_{x, \mu, \omega})}
\leq
  \int_0^t \int_0^1 \int_0^\infty \mu v_\omega |\phi_0(s, \mu, \omega)| \dmu \domega \ds.
\end{align*}

\Ni (b) Since the system is linear and $(f, g) = (m_0 \xi, \,\, \gamma_0 m_0 \xi)$ is a solution to the coupled system except the initial data, we only need to show that $f, g$ are non-positive given $\phi \leq 0$ and non-positive initial data. Let $f^+, g^+$ be the positive part of $f, g$ defined by
\begin{align*}
   f^+ = \max\{f, \, 0\}, 
\qquad
  g^+ = \max\{g, \, 0\}.
\end{align*}
Then they satisfy 
\begin{align} \label{sys:full-f-plus}
\begin{split}
 &  \del_t f^+ + \mu v_\omega \del_x f^+
\leq \frac{1}{\tau(\omega)} \vpran{\xi(\omega) \dvint{f^+} -  f^+}, 
\\
& f^+ |_{x=0} = 0, \hspace{4.9cm} \mu > 0,   
\\
&  f^+ |_{x=1} \leq \eta_1 f^+(1, -\mu, \omega) + \zeta_1 g^+(1, \mu, \omega), \qquad \mu < 0,  
\\
& f^+|_{t=0} \leq 0
\end{split}
\end{align}
and
\begin{align} \label{sys:full-g-plus}
\begin{split}
 &  \del_t g^+ + \mu v_\omega \del_x g^+ 
\leq \frac{1}{\tau(\omega)} \vpran{\xi(\omega) \dvint{g^+} -  g^+}, 
\\
& g^+ |_{x=1} \leq \eta_2 f^+(1, \mu, \omega) + \zeta_2 g^+(1, -\mu, \omega), \hspace{2.4cm} \mu > 0,   
\\
&  g^+ |_{x=L} \leq \alpha_0 \xi(\omega) \int_0^1 \int_0^\infty \mu \, v_\omega \, g^+(t, 1, \mu, \omega) \domega \dmu, \qquad \mu < 0,  
\\
& g^+|_{t=0} \leq 0.  
\end{split}
\end{align}
Then by the same argument as for the $L^1$ bound (replacing $\abs{f}, \abs{g}$ with $f^+, g^+$), we get
\begin{align*}
   \norm{f^+}_{L^\infty(0, t; L^1_{x, \mu, \omega})}
+ \norm{g^+}_{L^\infty(0, t; L^1_{x, \mu, \omega})}
\leq 0,
\end{align*}
which implies that $f, g \leq 0$. We thereby finish the proof by the comment at the beginning of (b). 

\smallskip
\Ni (c) Before applying the Riesz-Thorin Theorem using the $L^1-$ and $L^\infty-$bounds in parts (a) and (b), we normalize $(f, g)$ by 
\begin{align*}
    (\bar f, \ \bar g) = (f/\xi, \ g/\xi),
\end{align*}
such that the $L^1$ and $L^\infty$ spaces will have the same measure $\xi \dmu \domega$. 
Rewrite~\eqref{sys:full-f}-\eqref{sys:full-g} in terms of $(\bar f, \bar g)$:
\begin{align} \label{sys:full-f-scale}
\begin{split}
 &  \del_t \bar f + \mu v_\omega \del_x \bar f 
= \frac{1}{\tau(\omega)} \vpran{\dvint{\bar f}_\xi -  \bar f}, \hspace{1.9cm} x\in[0,1],
\\
& \bar f |_{x=0} = \phi/\xi, \hspace{5.3cm} \mu > 0,   
\\
&  \bar f|_{x=1} = \eta_1(\omega) \bar f(1, -\mu, \omega) + \zeta_1(\omega) \bar g(1, \mu, \omega), \qquad \mu < 0,  
\\
& \bar f|_{t=0} = 0
\end{split}
\end{align}
and
\begin{align} \label{sys:full-g-scale}
\begin{split}
 &  \del_t \bar g + \mu v_\omega \del_x \bar g 
= \frac{1}{\tau(\omega)} \vpran{\dvint{\bar g}_\xi -  \bar g} ,\hspace{2.9cm} x\in[1,L],
\\
& \bar g |_{x=1} = \eta_2(\omega) \bar f(1, \mu, \omega) + \zeta_2(\omega) \bar g(1, -\mu, \omega), \hspace{1.5cm} \mu > 0,   
\\
&  \bar g |_{x=L} = \alpha_0 \int_0^1 \int_0^\infty \mu \, v_\omega \, \bar g(t, L, \mu, \omega) \xi(\omega) \domega \dmu, \qquad \mu < 0,  
\\
& \bar g|_{t=0} = 0,
\end{split}
\end{align}
where $\dvint{h}_\xi = \int_{-1}^1 \int_0^\infty h \, \xi \domega \dmu$.
For a given $\phi$, if $(f, g)$ is the solution to~\eqref{sys:full-f}-\eqref{sys:full-g}, then $(\bar f, \bar g)$ is the solution to~\eqref{sys:full-f-scale}-\eqref{sys:full-g-scale}. Define the linear mapping 
\begin{align*}
   T (\phi/\xi) = (\bar f, \bar g).
\end{align*}
By the bounds in parts (a) and (b), for each $t > 0$, we have
\begin{align*}
  \norm{T (\phi / \xi)}_{L^1(\xi \dx \dmu \domega)}
 =  \norm{(\bar f, \ \bar g)}_{L^1(\xi \dx \dmu \domega)}
 \leq
   \norm{\phi/\xi}_{L^1(0, t; L^1(\xi \dmu\domega))}
\end{align*}
and
\begin{align*}
  \norm{T (\phi / \xi)}_{L^\infty(\xi \dx \dmu \domega)}
 =  \norm{(\bar f, \ \bar g)}_{L^\infty(\xi \dx \dmu \domega)}
 \leq
  (1 + \gamma_0) \norm{\phi/\xi}_{L^\infty(0, t; L^\infty(\xi \dmu\domega))}.
\end{align*}
By the Riesz-Thorin interpolation theorem, we derive that for any $p \in [1, \infty]$, the operator
\begin{align*}
  T: L^\infty(0, t; L^p(\xi \dx\dmu\domega))
  \to 
     L^p (0, t; L^p(\xi \dx\dmu\domega))
     \times L^p (0, t; L^p(\xi \dx\dmu\domega)) 
\end{align*}
is bounded with the bound
\begin{align*}
  \norm{(\bar f, \ \bar g)}_{L^p(\xi \dx \dmu \domega)}
  = \norm{T (\phi / \xi)}_{L^p(\xi \dx \dmu \domega)}
 \leq
  (1 + \gamma_0)^{1/p} \norm{\phi/\xi}_{L^p(0, t; L^p(\xi \dmu\domega))}, 
\quad 
  p \in [1, \infty],
\end{align*}
which when written in terms of $(f, g)$ gives the bound in~\eqref{bound:L-p}. 
\end{proof}

At the end of this section we recall some solution formulae for later use. Let $f$ be the solution to the~system
\begin{align} \label{sys:general-f}
\begin{split}
 &  \del_t f + \mu v_\omega \del_x f 
= -\frac{1}{\tau(\omega)} f + \frac{\xi(\omega)}{\tau(\omega)} h(t, x),
\\
& f |_{x=0} = \phi(t, \mu, \omega), \hspace{3.3cm} \mu > 0,   
\\
&  f|_{x=1} = \eta_1 f(1, -\mu, \omega) + \zeta_1 g(1, \mu, \omega), \qquad \mu < 0,  
\\
& f|_{t=0} = 0,
\end{split}
\end{align}
then $f$ can be explicitly solved via the method of characteristics. In particular, the incoming part of $f$ at $x=1$ is
\begin{align} \label{soln:x-incom-1_appendix}
  f(t, 1, \mu, \omega) \, \One_{\mu < 0} 
= \begin{cases}
   \zeta_1 g(t, 1, \mu, \omega)
   + \eta_1 \int_0^t e^{-\frac{y}{\tau(\omega)}}
   \frac{\xi(\omega)}{\tau(\omega)}
   h(t-y, 1 + \mu v_\omega y) \dy, 
   & 0 < t < \frac{1}{|\mu| v_\omega}, \ \mu < 0,  \\[10pt]
   \zeta_1 g(t, 1, \mu, \omega)
   + \eta_1 \int_0^{\frac{1}{|\mu| v_\omega}} e^{-\frac{y}{\tau(\omega)}}
   \frac{\xi(\omega)}{\tau(\omega)}
   h(t-y, 1 + \mu v_\omega y) \dy \\
   \hspace{3.9cm}+ \eta_1 \phi\vpran{t + \frac{1}{\mu \omega}, - \mu, \omega} e^{-\frac{1}{\tau(\omega)} \frac{1}{|\mu| v_\omega}}, 
   & t > \frac{1}{|\mu| v_\omega}, \quad \mu < 0.
    \end{cases}
\end{align}
Solving backwards in $x$ gives the outgoing data at $x=0$ as
\begin{align} \label{soln:x-out-0_appendix}
 f(t, 0, \mu, \omega) \, \One_{\mu < 0} 
= \begin{cases}
   \int_0^t e^{-\frac{1}{\tau(\omega)} y}
   \frac{\xi(\omega)}{\tau(\omega)}
   h(t-y, -\mu v_\omega y) \dy, 
   & 0 < t < \frac{1}{|\mu| v_\omega},  \\[10pt]
   \int_0^{\frac{1}{|\mu| v_\omega}} e^{-\frac{1}{\tau(\omega)} y}
   \frac{\xi(\omega)}{\tau(\omega)}
   h(t-y, -\mu v_\omega y) \dy
   + f\vpran{t + \frac{1}{\mu \omega}, 1, \mu, \omega}
     e^{-\frac{1}{\tau(\omega)} \frac{1}{|\mu| v_\omega}}, 
   & t > \frac{1}{|\mu| v_\omega}. 
    \end{cases}
\end{align}

\section{Reconstruction theory}\label{sec:theory}
In this section we present the main part of the paper: reconstruction of the heat-reflection indices. 
As explained earlier, we take the singular decomposition approach by identifying the terms that contribute to the measurements at the leading order and showing the remainder terms are indeed negligible in the measurements. 

\subsection{Reconstruction strategy and results}\label{sec:strategy}
To reconstruct the reflection coefficient $\eta_1(\omega)$, we inject a rather singular source at the boundary concentrated at $\omega_0$ and $\mu_0$. The leading order of the solution with such a source will be an X-ray propagation concentrating at $(\mu_0,\omega_0)$. 
Due to the concentration, the leading-order energy will be reflected back to the physical boundary at a pre-designed time $\frac{2}{\mu_0 v_{\omega_0}}$. The measurement is thus chosen to be taken at that particular time to capture the leading order.



To be more specific, we decompose $f$ in equation~\eqref{sys:full-f} as
\begin{align}\label{eqn:f_decompose}
f=f_0+f_1\,,
\end{align}
with each component satisfying
\begin{align} \label{sys:f-0}
\begin{split}
 &  \del_t f_0 + \mu v_\omega \del_x f_0 
= - \frac{1}{\tau(\omega)} f_0, 
\\
& f_0 |_{x=0} = \phi(t, \mu, \omega), \hspace{2.1cm} \mu > 0,   
\\
&  f_0|_{x=1} = \eta_1 f_0(1, -\mu, \omega), \hspace{1.4cm} \mu < 0,  
\\
& f_0 |_{t=0} = 0
\end{split}
\end{align}
and
\begin{align} \label{sys:f-1}
\begin{split}
 &  \del_t f_1 + \mu v_\omega \del_x f_1
= - \frac{1}{\tau(\omega)} f_1
   + \frac{\xi(\omega)}{\tau(\omega)} \dvint{f}, 
\\
& f_1 |_{x=0} = 0, \hspace{5.4cm} \mu > 0,   
\\
&  f_1 |_{x=1} = \eta_1 f_1(1, -\mu, \omega)
                         + \zeta_1 g(1, \mu, \omega), \hspace{1.4cm} \mu < 0,  
\\
& f_1 |_{t=0} = 0. 
\end{split}
\end{align}
Note that $f_0$ takes all the information of the incoming data. While $f_1$ has the trivial boundary condition at $x=0$, and it picks up all the scattering involving $\dvint{f}$ as a source term. 
Choose the incoming data $\phi$ as
\begin{align}\label{eqn:phi}
   \phi(t, \mu, \omega)
= \tfrac{1}{\Eps^3} \phi_0 \vpran{\tfrac{t}{\Eps}}
    \phi_0 \vpran{\tfrac{\mu - \mu_0}{\Eps}}
    \phi_0 \vpran{\tfrac{\omega - \omega_0}{\Eps}}, 
\end{align}
where $\phi_0$ is a smooth cutoff function satisfying 
\begin{align} \label{def:phi-0}
   \int_{\R} \phi_0(z) \dz = 1,
\qquad
   \phi_0(z)
= \begin{cases}
     0, & \text{if $z < 0$ or $z \geq 1$}, \\[2pt]
     0 < \phi(0) < 1, & \text{if $z \in (0, 1)$}, \\[2pt]
     \text{smooth}, & \text{for $z \in (0, 1)$}.
   \end{cases}
\end{align}
Here $\mu_0$, $\omega_0$ are fixed and $\Eps$ is set to be small to ensure the concentration. Conditions on these parameters will be specified later. Define the measuring operator as follows: for any $h = h(t, \mu, \omega)$, let
\begin{align}\label{eqn:calM}
   \CalM(h)
= \int_0^\infty \psi_0 \vpran{\frac{t- t_1}{\Eps}}
   \vpran{\int_0^\infty \int_{-1}^0 h(t, \mu, \omega) \dmu \domega} \dt, 
\end{align}
where $t_1, \psi_0$ satisfies 
\begin{align} \label{def:psi-0}
  t_1 = \tfrac{2}{\mu_0 v_{\omega_0}}, 
\qquad
   \int_{\R} \psi_0(z) \dz = 1,
\qquad
   \psi_0(z)
= \begin{cases}
     0, & \text{if $z < -1$ or $z \geq 1$}, \\[2pt]
     0 < \phi(0) < 1, & \text{if $z \in (-1, 1)$}, \\[2pt]
     \text{smooth}, & \text{for $z \in \R$}. 
   \end{cases}
\end{align}
This particular way of choosing $\psi_0$ is to specifically single out the information at time $t_1$. This is the time when the concentrated boundary data injected into $f_0$ gets bounced back by the interface and returns to $x=0$. 

The main assumptions we make for the reconstruction are
\begin{itemize}
\item[(A1)] The group velocity $v_\omega$ is differentiable and bounded: there exists $v_0 > 0$ such that
\begin{align} \label{assump:velo}
   v_0 : = \sup_{\omega} v_\omega < \infty. 
\end{align}

\item[(A2)] There exists $p_0 \in (1, 3/2)$ such that the group velocity $v_\omega$ satisfies
\begin{align}\label{assump:p0}
   \int_0^\infty \frac{\xi(\omega)}{\tau(\omega)} \frac{1}{v_\omega^{1/p_0}} \domega < \infty. 
\end{align}
Denote $p'_0$ as the H\"older conjugate of $p_0$.

\item[(A3)] The relaxation time $\tau(\omega)$ is bounded away from zero: there exists $\tau_0 > 0$ such that
\begin{align} \label{assump:tau}
    \tau(\omega) \geq \tau_0 > 0.
\end{align}

\item[(A4)] The width of the right region is large enough:
\begin{align} \label{assump:L}
    L \geq \frac{v_0}{\mu_0 v_{\omega_0}} + \frac{v_0}{2} + 1.
\end{align}
\end{itemize}
The last assumption (A4) can be viewed either as a constraint for the range of $\omega_0$ if $L$ is fixed, or it can be viewed as a condition for choosing $L$ for a specific range of $\omega_0$. This assumption is to guarantee there is a sufficient delay in time for majority of the energy back from $x=L$ to be sensed at the surface $x=0$. If the boundary condition at $x=L$ is absorbing-type, meaning all the heat hitting the boundary $L$ will simply be released from the system to the outer space, the assumption can be removed.

The main result of this paper is
\begin{thm}\label{thm:main}
Let assumptions (A1)-(A4) hold true and choose $\phi_0, \psi_0$ as described in~\eqref{def:phi-0} and \eqref{def:psi-0}. Let $(f, g)$ be the solution to the system~\eqref{sys:full-f}-\eqref{sys:full-g} with incoming data $\phi_0$. Then
\begin{align*}
    \lim_{\Eps \to 0} \CalM(f) 
= \eta_1(\omega_0) \, C_{\phi_0, \psi_0}, 
\end{align*}
where the constant $C_{\phi_0, \psi_0}$ can be pre-calculated as \begin{align}\label{eqn:C_0}
\begin{split}
   C_{\phi_0, \psi_0}
&= e^{- \frac{1}{\tau(\omega_0)} \frac{2}{\mu_0 v_{\omega_0}}}
  \int_0^\infty \int_0^1 \int_0^\infty
  \psi_0 \vpran{t + \frac{2}{(\mu_0 v_{\omega_0})^2}
   \vpran{v_{\omega_0} \mu + \mu_0 v'(\omega_0) \omega}}
\\
& \hspace{6cm} \times
   \phi_0(t) \phi_0(\mu) \phi_0(\omega) \dt \dmu \domega.  \qedhere
\end{split}
\end{align}
This gives the reconstruction of $\eta_1(\omega_0)$ as
\begin{align*}
   \eta_1(\omega_0) 
= \frac{1}{C_{\phi_0, \psi_0}} \lim_{\Eps \to 0} \CalM(f).
\end{align*}
The recovery of $(\eta_2, \zeta_1, \zeta_2)$ follows from~\eqref{assump:coeff}.
\end{thm}

The proof of Theorem~\ref{thm:main} will be divided into two steps. Noting that $f=f_0+f_1$ in~\eqref{eqn:f_decompose}, we decompose the measurement accordingly as
\begin{align*}
\CalM(f) = \CalM(f_0) + \CalM(f_1)\,.
\end{align*}
We will show below that as $\epsilon\to0$, $\CalM(f) = \CalM(f_0)$ and $\CalM(f_1)= 0$, as presented in Propositions~\ref{prop:f0} and~\ref{prop:f1} respectively. 

\begin{prop}\label{prop:f0}
Let assumptions (A1)-(A4) hold true and choose $\phi_0, \psi_0$ as described in~\eqref{def:phi-0} and \eqref{def:psi-0}. Let $f_0$ be the solution to~\eqref{sys:f-0}. Then
\begin{align} \label{eqn:f_0_limit}
    \lim_{\Eps \to 0} \CalM(f_0) 
= \eta_1(\omega_0) C_{\phi_0, \psi_0}\,.
\end{align}
\end{prop}
\begin{proof}

Directly solving~\eqref{sys:f-0} along its trajectories, we obtain the outgoing data at $x=0$ as
\begin{align*}
   f_0(t, 0, \mu, \omega) \One_{\mu < 0} 
= \begin{cases}
   0, 
   & 0 < t < \frac{2}{|\mu| v_\omega},  \\[10pt]
   \eta_1 \phi\vpran{t + \frac{2}{\mu v_\omega}, -\mu, \omega}
     e^{-\frac{2}{\tau(\omega)} \frac{1}{|\mu| v_\omega}} \One_{\mu < 0}, 
   & t > \frac{2}{|\mu| v_\omega}.
    \end{cases}
\end{align*}
Therefore, the contribution of $f_0$ toward the measurement is 
\begin{align} \label{meas:f-0}
\begin{split}
  \CalM(f_0) 
&= \int_0^\infty \int_{-1}^0 \int_{\frac{2}{|\mu| \omega}}^\infty
  \eta_1(\omega) \psi_0 \vpran{\frac{t - t_1}{\Eps}}
   e^{-\frac{1}{\tau_\omega} \frac{2}{|\mu| v_\omega}}
    \phi \vpran{t + \tfrac{2}{\mu v_\omega}, -\mu, \omega} \dt \dmu \domega  
\\
& \quad \hspace{-1cm}= \frac{1}{\Eps^3} \int_0^\infty \int_0^1 \int_{\frac{2}{|\mu| \omega}}^\infty
   \eta_1(\omega) \psi_0 \vpran{\frac{t - t_1}{\Eps}}
   e^{-\frac{1}{\tau_\omega} \frac{2}{\mu v_\omega}}
   \phi_0 \vpran{\frac{t - \tfrac{2}{\mu v_\omega}}{\Eps}}
   \phi_0 \vpran{\frac{\mu - \mu_0}{\Eps}}
   \phi_0 \vpran{\frac{\omega - \omega_0}{\Eps}} \dt \dmu \domega. 
\end{split}
\end{align}
To compute the limit, we apply the Lebesgue Dominated Convergence Theorem (LDCT). First we make the change of variables by letting
\begin{align*}
   \tilde t = \frac{t - \frac{2}{\mu v_\omega}}{\Eps},
\qquad
  \tilde \mu = \frac{\mu - \mu_0}{\Eps},
\qquad
  \tilde \omega = \frac{\omega - \omega_0}{\Eps}. 
\end{align*}
Applying these new variables in the integral in ~\eqref{meas:f-0} gives, and utilizing the support information of $\phi$, we have:
\begin{align*}
  \CalM(f_0)
= \int_0^1 \int_0^1 \int_0^1
    \eta_1(\omega_0 + \Eps \, \tilde \omega)
    e^{- \frac{1}{\tau(\omega)} \frac{2}{\mu v_\omega}}
    \psi_0\vpran{\tilde t  + \frac{\frac{2}{\mu_0 v_{\omega_0}} \vpran{1 - \frac{\mu_0 v_{\omega_0}}{\mu v_\omega}} }{\Eps}} \,
    \phi_0(\tilde t) \, \phi_0 (\tilde \mu) \, \phi_0 (\tilde \omega) 
    {\, \rm d} \tilde t {\, \rm d} \tilde \mu {\,\rm d} \tilde \omega,
\end{align*}
where for the simplicity of notations, we have kept $(\mu, \omega)$ as a function of $(\tilde \mu, \, \tilde \omega)$ in the exponential term. By the continuity of the functions, we obtain the pointwise limits as
\begin{align*}
   \lim_{\Eps \to 0} \, \eta_1(\omega_0 + \Eps \, \tilde \omega) \,
    e^{- \frac{1}{\tau(\omega)} \frac{2}{\mu v_\omega}}
 = \eta_1 (\omega_0) \, e^{- \frac{1}{\tau(\omega_0)} \frac{2}{\mu_0 v_{\omega_0}}}
\end{align*}
and
\begin{align*}
  \lim_{\Eps \to 0} \frac{\frac{2}{\mu_0 v_{\omega_0}} \vpran{1 - \frac{\mu_0 v_{\omega_0}}{\mu v_\omega}} }{\Eps}
= \lim_{\Eps \to 0} \frac{\frac{2}{\mu_0 v_{\omega_0}} \vpran{1 - \frac{\mu_0 v_{\omega_0}}{(\mu_0 + \Eps \tilde \mu) v(\omega_0 + \Eps \tilde \omega)}} }{\Eps}
= \frac{2}{(\mu_0 v_{\omega_0})^2}
   \vpran{v_{\omega_0} \tilde \mu + \mu_0 v'(\omega_0) \tilde \omega}.
\end{align*}
where $v'$ denotes the derivative of $v$. By the boundedness of $\eta_1, \psi_0, \phi_0$ and the boundedness of the integration domain, the LDCT applies and gives
\begin{align*}
  \lim_{\Eps \to 0} \CalM(f_0)
&= \eta_1 (\omega_0) \, e^{- \frac{1}{\tau(\omega_0)} \frac{2}{\mu_0 v_{\omega_0}}}
  \int_0^\infty \int_0^1 \int_0^\infty
  \psi_0 \vpran{t + \frac{2}{(\mu_0 v_{\omega_0})^2}
   \vpran{v_{\omega_0} \mu + \mu_0 v'(\omega_0) \omega}}
\\
& \hspace{5cm} \times
   \phi_0(t) \phi_0(\mu) \phi_0(\omega) \dt \dmu \domega,
\end{align*}
where we have dropped the tilde signs and extended the integral domain to the full one using the supports of $\psi_0$ and $\phi_0$. This verifies the formula for $C_{\phi_0, \psi_0}$ in~\eqref{eqn:C_0} and proves~\eqref{eqn:f_0_limit}.
\end{proof}

We are left to show that $\CalM(f_1)$ vanishes in the limit. Recall that $f_1$ satisfies~\eqref{sys:f-1}. Using formula~\eqref{soln:x-out-0_appendix}, we obtain that for $x=0$ and $\mu < 0$,
\begin{align*}
   f_1(t, 0, \mu, \omega) \One_{\mu < 0}
= \begin{cases}
   \int_0^t e^{-\frac{y}{\tau(\omega)}}
   \frac{\xi(\omega)}{\tau(\omega)}
   \dvint{f}(t-y, -\mu v_\omega y) \dy, 
   & 0 < t < \frac{1}{|\mu| v_\omega},  \\[10pt]
   \int_0^{\frac{1}{|\mu| v_\omega}} e^{-\frac{y}{\tau(\omega)}}
   \frac{\xi(\omega)}{\tau(\omega)}
   \dvint{f}(t-y, -\mu v_\omega y) \dy
   + f_1 \vpran{t + \frac{1}{\mu v_\omega}, 1, \mu, \omega}
     e^{-\frac{1}{\tau(\omega)} \frac{1}{|\mu| v_\omega}}, 
   & t > \frac{1}{|\mu| v_\omega},
    \end{cases}
\end{align*}
where by~\eqref{soln:x-incom-1_appendix} the incoming part at $x=1$ is 
\begin{align} \label{soln:x-incom-1}
  f_1(t, 1, \mu, \omega) \, \One_{\mu < 0} 
= \begin{cases}
   \zeta_1 g(t, 1, \mu, \omega)
   + \eta_1 \int_0^t e^{-\frac{y}{\tau(\omega)}}
   \frac{\xi(\omega)}{\tau(\omega)}
   \dvint{f}(t-y, 1 + \mu v_\omega y) \dy, 
   & 0 < t < \frac{1}{|\mu| v_\omega},  \\[10pt]
   \zeta_1 g(t, 1, \mu, \omega)
   + \eta_1 \int_0^{\frac{1}{|\mu| v_\omega}} e^{-\frac{y}{\tau(\omega)}}
   \frac{\xi(\omega)}{\tau(\omega)}
   \dvint{f}(t-y, 1 + \mu v_\omega y) \dy, 
   & t > \frac{1}{|\mu| v_\omega}\,.
    \end{cases}
\end{align}
More explicitly, if $0 < t < \tfrac{1}{|\mu| v_\omega}$, then
\begin{align} \label{soln:f-1-1}
   f_1(t, 0, \mu, \omega) \One_{\mu < 0}
= \int_0^t e^{-\frac{y}{\tau(\omega)}}
   \frac{\xi(\omega)}{\tau(\omega)}
   \dvint{f}(t-y, \, -\mu v_\omega y) \dy\,.
\end{align}
If $\tfrac{1}{|\mu| v_\omega} < t < \tfrac{2}{|\mu| v_\omega}$, then, calling the first item in~\eqref{soln:x-incom-1}, we have:
\begin{align} \label{soln:f-1-2}
\begin{split}
   f_1(t, 0, \mu, \omega) \One_{\mu < 0}
&= \int_0^{\frac{1}{|\mu| v_\omega}} e^{-\frac{y}{\tau(\omega)}}
   \frac{\xi(\omega)}{\tau(\omega)}
   \dvint{f}(t-y, -\mu v_\omega y) \dy
   + \zeta_1 g(t + \tfrac{1}{\mu v_\omega}, 1, \mu, \omega) 
      e^{-\frac{1}{\tau(\omega)} \frac{1}{|\mu| v_\omega}}
\\
& \quad \,
  + \eta_1 e^{-\frac{1}{\tau(\omega)} \frac{1}{|\mu| v_\omega}} 
     \int_0^{t+\frac{1}{\mu v_\omega}} e^{-\frac{y}{\tau(\omega)}}
   \frac{\xi(\omega)}{\tau(\omega)}
   \dvint{f}(t + \tfrac{1}{\mu v_\omega} -y, \, 1 + \mu v_\omega y) \dy. 
\\
& = : f_{1,1} + f_{1,2} + f_{1,3}\,.
\end{split}
\end{align}
Finally, if $t > \frac{2}{|\mu| v_\omega}$, then, calling the second item in~\eqref{soln:x-incom-1}, we have:
\begin{align} \label{soln:f-1-3}
\begin{split}
   f_1(t, 0,  \mu, \omega) \One_{\mu < 0}
& = \int_0^{\frac{1}{|\mu| v_\omega}} e^{-\frac{y}{\tau(\omega)} }
   \frac{\xi(\omega)}{\tau(\omega)}
   \dvint{f}(t-y, -\mu v_\omega y) \dy
   + \zeta_1 g(t + \tfrac{1}{\mu v_\omega}, 1, \mu, \omega) 
      e^{-\frac{1}{\tau(\omega)} \frac{1}{|\mu| v_\omega}}
\\
& \quad \,
   + \eta_1 e^{-\frac{1}{\tau(\omega)} \frac{1}{|\mu| v_\omega}}
      \int_0^{\frac{1}{|\mu| v_\omega}} e^{-\frac{y}{\tau(\omega)}}
   \frac{\xi(\omega)}{\tau(\omega)}
   \dvint{f}(t + \tfrac{1}{\mu v_\omega} -y, \, 1 + \mu v_\omega y) \dy
\\
& = : f_{1,4} + f_{1,5} + f_{1,6}\,.
\end{split}
\end{align}
Note the only difference between $f_{1,3}$ and $f_{1,6}$ lies in the lower and upper bound of the integration.

Accordingly, the contribution of $f_1$ to the measurement can be broken up as
\begin{align} \label{decomp:M-f-1}
\begin{split}
   \CalM(f_1)
&= \int_0^\infty \int_{-1}^0 \int_{0}^\infty
   \psi_0 \vpran{\frac{t - t_1}{\Eps}}
   f_1(t, 0, \mu, \omega) \dt \dmu \domega.
\\
& = \CalM \vpran{f_1 \One_{0 < t < \tfrac{1}{|\mu| v_\omega}}}
      + \underbrace{\CalM \vpran{f_1 \One_{\tfrac{1}{|\mu| v_\omega} < t < \tfrac{2}{|\mu| v_\omega}}}}_{\CalM(f_{1,1})+\CalM(f_{1,2})+\CalM(f_{1,3})}
      + \underbrace{\CalM \vpran{f_{1} \One_{t > \tfrac{2}{|\mu| v_\omega}}}}_{\CalM(f_{1,4})+\CalM(f_{1,5})+\CalM(f_{1,6})}\,.
\end{split}
\end{align}
Among these terms, $\CalM(f_{1,2})$ and $\CalM(f_{1,5})$ contain the contributions from $g$. $\CalM(f_{1,1})$, $\CalM(f_{1,3})$, $\CalM(f_{1,4})$ and $\CalM(f_{1,6})$ depend on $\dvint{f}$. In Propositions~\ref{prop:f1}, we quantify the smallness of each one of them. The proofs of these terms follow a similar pattern shown in
\begin{prop} \label{prop:general}
Suppose $\Eps < 1$ and assumptions (A1)-(A4) hold. Suppose $(f, g)$ is the solution to~\eqref{sys:full-f}-\eqref{sys:full-g} with incoming data $\phi_0$ defined in~\eqref{def:phi-0}.
Then there exists a constant $C$ depending on $p_0$, $\phi_0$, $\gamma_0$ and $t_1$ such that
 \begin{align} \label{bound:general-f}
   \abs{\int_0^\infty \int_{-1}^0 \int_0^\infty \int_0^t
   \psi_0 \vpran{\frac{t - t_1}{\Eps}} \frac{\xi(\omega)}{\tau(\omega)}
   \dvint{f}(t - y, \, -\mu v_\omega y) \One_{|\mu| v_\omega y < 1} \dy \dt}
\leq  
  C \Eps^{1-\frac{3}{p'_0}} \to 0
\qquad
  \text{as $\Eps \to 0$.}
\end{align}
where $p'_0$ is the H\"{o}lder conjugate of the constant $p_0$ in assumption (A2). Similar bound holds for $g$. 
\end{prop}
\begin{proof}
First, by applying~\eqref{bound:L-p} in Proposition~\ref{prop:f-L-p} to system~\eqref{sys:full-f}-\eqref{sys:full-g}, we have 
\begin{align} \label{bound:f-L-p}
   \norm{f(t, \cdot, \cdot, \cdot)}_{L^p(\xi^{1-p} \dx \dmu \domega)}
\leq
  C_{p, \gamma_0} \norm{\phi}_{L^p(0, t; L^p(\xi^{1-p} \dmu\domega))}
\leq
  C_{p, \phi_0, \gamma_0} \Eps^{-3(1 - \tfrac{1}{p})}, 
\qquad
  \forall \, p \in [1, \infty)\,.
\end{align}
By the support of $\psi_0$, we have
\begin{align} \label{bound:general-1}
\begin{split}
& \quad \,
   \abs{\int_0^\infty \int_{-1}^0 \int_{0}^\infty \int_0^t
   \psi_0 \vpran{\frac{t - t_1}{\Eps}}
   \frac{\xi(\omega)}{\tau(\omega)}
   \dvint{f}(t-y, -\mu v_\omega y)  \One_{|\mu| v_\omega y < 1}
   \dy \dt \dmu \domega}
\\
& \leq 
  \int_0^\infty \frac{\xi(\omega)}{\tau(\omega)} 
  \vpran{\underbrace{\int^{t_1+\Eps}_{t_1-\Eps} \int_0^t \int_0^1 
   \dvint{|f|}(t-y, \mu v_\omega y) \One_{|\mu| v_\omega y < 1}
   \dmu \dy \dt}_{\mathrm{I}}} \domega\,.
\end{split}
\end{align}
The key step is to make the change of variables $\mu \to x$ by 
\begin{align*}
   x = \mu v_\omega y. 
\end{align*}
Using the new variable $x$, we have
\begin{align} \label{bound:general-2}
\mathrm{I} = \int^{t_1+\Eps}_{t_1-\Eps} \int_0^t  \frac{1}{v_\omega y}
  \vpran{\int_0^{v_\omega y}
   \dvint{|f|}(t-y, x) \One_{x < 1}
   \dx} \dy \dt. 
\end{align}
To estimate the $\dx$-integral and reduce the singularity of $y$ at zero, we use the $L^p$-bound of $f$. By H\"{o}lder's inequality and~\eqref{bound:f-L-p}:
\begin{align} \label{bound:L-p-step}
& \quad \,
   \int_0^{v_\omega y}
   \dvint{|f|}(t-y, x) \One_{x < 1} \dx
= \int_0^{v_\omega y} \int_{-1}^1 \int_0^\infty
      \abs{f (t-y, x, \mu', \omega')} \One_{x < 1}
      \domega' \dmu' \dx \nn
\\
& \leq
   \vpran{\int_0^1 \int_{-1}^1 \int_0^\infty
      \abs{f (t-y, x, \mu', \omega')}^p
      \xi^{1-p}(\omega') \domega' \dmu' \dx}^{\frac{1}{p}}
   \vpran{\int_0^{v_\omega y} \int_{-1}^1 \int_0^\infty
      \xi(\omega') \domega' \dmu' \dx}^{\frac{1}{p'}} \nn
\\
& \leq  
   C_p \vpran{\sup_{t \in [0, t_1 + \Eps]}\norm{f(t, \cdot, \cdot, \cdot)}_{L^p(\xi^{1-p} \domega \dmu \dx)}} \vpran{v_\omega y}^{1/p'}
\leq
  C_{p, \phi_0, \gamma_0} \Eps^{-\frac{3}{p'}}
  \vpran{v_\omega y}^{1/p'},
\qquad
  \forall p \in (1, \infty).
\end{align}
The last inequality calls for the $L^p$ estimate of $f$ in~\eqref{bound:f-L-p}. Inserting this estimate into~\eqref{bound:general-1} gives
\begin{align*}
& \quad \,
   \abs{\int_0^\infty \int_{-1}^0 \int_{0}^\infty \int_0^t
   \psi_0 \vpran{\frac{t - t_1}{\Eps}}
   \frac{\xi(\omega)}{\tau(\omega)}
   \dvint{f}(t-y, -\mu v_\omega y)  \One_{|\mu| v_\omega y < 1}
   \dy \dt \dmu \domega}
\\
&\leq
  C_{p, \phi_0, \gamma_0} \Eps^{-\frac{3}{p'}} \int_0^\infty \frac{\xi(\omega)}{\tau(\omega)} \frac{1}{v_\omega^{1/p}}
  \vpran{\int_{t_1-\Eps}^{t_1+\Eps} \int_0^t \frac{1}{y^{1/p}} \dy \dt} \domega
\\
&\leq
  C_{p, \phi_0, \gamma_0, t_1} \Eps^{1-\frac{3}{p'}}
  \vpran{\int_0^\infty \frac{\xi(\omega)}{\tau(\omega)} \frac{1}{v_\omega^{1/p}} \domega}. 
\end{align*}
Hence if we choose $p = p_0 < 3/2$ as prescribed in assumption (A2), then $1-\frac{3}{p'} > 0$ and~\eqref{bound:general-f} holds. Since the only ingredient from $f$ is its $L^p$-bound and $g$ satisfies a similar bound, ~\eqref{bound:general-f} holds with $f$ replaced by $g$. 
\end{proof}

Now we apply Proposition~\ref{prop:general} to show the smallness of each term in $\CalM(f_1)$.

\begin{prop}\label{prop:f1}
Suppose $\Eps < 1$ and $t_1 = \frac{2}{\mu_0 v_{\omega_0}}$. Let $p_0 \in (1, \frac{3}{2})$ be the constant in assumption (A2) and let $p'_0$ be its H\"{older} conjugate. Then there exists a constant $C$ depending on $p_0$, $\phi_0$, $\gamma_0$ and $t_1$ such that
 \begin{align} \label{bound:M-f-1-total}
 \CalM \vpran{f_1} 
\leq  
  C \Eps^{1-\frac{3}{p'_0}} \to 0
\qquad
  \text{as $\Eps \to 0$.}
\end{align}
\end{prop}
\begin{proof}
We will treat each term in~\eqref{soln:f-1-1}-\eqref{soln:f-1-3} separately. Start with the terms related to $\dvint{f}$. We show that they all fit into the structure in Proposition~\ref{prop:general}. First,
\begin{align} \label{est:M-f-1-1}
\begin{split}
  \CalM \vpran{f_1 \One_{0 < t < \tfrac{1}{|\mu| v_\omega}}}
& = \int_0^\infty \int_{-1}^0 \int_{0}^{\infty}
   \psi_0 \vpran{\frac{t - t_1}{\Eps}}
   f_1(t, 0, \mu, \omega) \One_{t < \frac{1}{|\mu| v_\omega}}\dt \dmu \domega 
\\
& = \int_0^\infty \int_0^1 \int_{0}^{\infty} \int_0^t
   \psi_0 \vpran{\frac{t - t_1}{\Eps}}
   \frac{\xi(\omega)}{\tau(\omega)}
   \dvint{f}(t-y, \mu v_\omega y)  \One_{t < \frac{1}{|\mu| v_\omega}} \dy \dt \dmu \domega,
\end{split}
\end{align}
which has the structure as in~\eqref{bound:general-f}. Hence Proposition~\ref{prop:general} applies. Next, 
\begin{align} \label{est:M-f-1-1-intm}
\begin{split}
  \CalM(f_{1,1})
&\leq
  \int_0^\infty \int_{-1}^0 \int_0^\infty 
  \psi_0\vpran{\tfrac{t-t_1}{\Eps}}
  \abs{f_{1,1}(t, 0, \mu, \omega)}
  \One_{\frac{1}{|\mu| v_\omega} < t < \frac{2}{|\mu| v_\omega}} \dt \dmu \domega
\\
&\leq
  \int_0^\infty \int_{-1}^0 \int_0^\infty  
  \int_0^{\frac{1}{|\mu| v_\omega}}
  \frac{\xi(\omega)}{\tau(\omega)}
  \psi_0\vpran{\tfrac{t-t_1}{\Eps}}
  \dvint{|f|}(t-y, |\mu| v_\omega y)
  \One_{\frac{1}{|\mu| v_\omega} < t < \frac{2}{|\mu| v_\omega}}
   \dy \dt \dmu \domega
\\
& = \int_0^\infty \int_{-1}^0 \int_0^\infty  
  \int_0^{t}
  \frac{\xi(\omega)}{\tau(\omega)}
  \psi_0\vpran{\tfrac{t-t_1}{\Eps}}
  \dvint{|f|}(t-y, |\mu| v_\omega y)
  \One_{\frac{1}{|\mu| v_\omega} < t < \frac{2}{|\mu| v_\omega}}
  \One_{0 < y < \frac{1}{|\mu| v_\omega}}
   \dy \dt \dmu \domega,
\end{split}
\end{align}
to which Proposition~\ref{prop:general} applies. To transform $\CalM(f_{1,3})$ into a form so that Proposition~\ref{prop:general} applies, we use the change of variables $y \to \tilde y$ where
\begin{align*}
    \tilde y = y - \frac{1}{\mu v_\omega}.
\end{align*}
 Then
\begin{align*}
  \CalM(f_{1,3})
& \leq
  \int_0^\infty \int_{-1}^0 \int_0^\infty
  \int_0^{t+\frac{1}{\mu v_\omega}}
  \frac{\xi(\omega)}{\tau(\omega)}
  \psi_0\vpran{\tfrac{t-t_1}{\Eps}}
  \dvint{|f|}(t+\tfrac{1}{\mu v_\omega} - y, 1+\mu v_\omega y) 
  \One_{\frac{1}{|\mu| v_\omega} < t < \frac{2}{|\mu| v_\omega}}\dy\dt\dmu\domega
\\
& = \int_0^\infty \int_{-1}^0 \int_0^\infty
  \int_{\frac{1}{|\mu| \omega}}^{t}
  \frac{\xi(\omega)}{\tau(\omega)}
  \psi_0\vpran{\tfrac{t-t_1}{\Eps}}
  \dvint{|f|}(t - \tilde y, 2+\mu v_\omega \tilde y) 
  \One_{\frac{1}{|\mu| v_\omega} < t < \frac{2}{|\mu| v_\omega}}
    {\, \rm d}\tilde y \dt\dmu\domega
\\
& = \int_0^\infty \int_{-1}^0 \int_0^\infty
  \int_0^t
  \frac{\xi(\omega)}{\tau(\omega)}
  \psi_0\vpran{\tfrac{t-t_1}{\Eps}}
  \dvint{|f|}(t - \tilde y, 2+\mu v_\omega \tilde y) 
  \One_{\frac{1}{|\mu| v_\omega} < t < \frac{2}{|\mu| v_\omega}}
  \One_{\tilde y > \frac{1}{|\mu| v_\omega}}
  {\, \rm d}\tilde y \dt\dmu\domega
\end{align*}
Make a further change of variables $\mu \to \tilde \mu$ by
\begin{align*}
  \tilde \mu 
= \mu \vpran{1 + \frac{2}{\mu v_\omega \tilde y}}
= \mu + \frac{2}{v_\omega \tilde y} = \mu\frac{\mu v_\omega y +1}{\mu v_\omega y -1}\,. 
\end{align*}
Since $\mu<0$, we have $|\tilde\mu|<|\mu|\leq 1$. Moreover, using the bounds $0<y<t+\frac{1}{\mu. v_\omega}$ and $t<\frac{2}{|\mu|v_\omega}$, we derive that $\mu v_\omega y + 1>0$ which gives $\tilde\mu>0$.
Hence,
\begin{align*}
   \CalM(f_{1,3})
\leq
  \int_0^\infty \int_0^\infty
  \int_0^t \int_{0}^1
  \frac{\xi(\omega)}{\tau(\omega)}
  \psi_0\vpran{\tfrac{t-t_1}{\Eps}}
  \dvint{|f|}(t - \tilde y, \,\, \tilde \mu v_\omega \tilde y) 
  {\,\rm d}\tilde \mu {\, \rm d} \tilde y\dt\domega,
\end{align*}
to which Proposition~\ref{prop:general} applies. Estimates for $\CalM(f_{1,4})$ and $\CalM(f_{1,6})$ are similar to $\CalM(f_{1,1})$ and $\CalM(f_{1,3})$ and the details are omitted to avoid repetition. Overall, we have
\begin{align} \label{bound:M-f-1-part-1}
   \CalM(f_{1, 1}) + \CalM(f_{1, 3}) + \CalM(f_{1, 4}) + \CalM(f_{1, 6})
\leq 
   C_{p, \phi_0, \gamma_0, t_1} \Eps^{1-\frac{3}{p'_0}}. 
\end{align}

We are left to estimate $\CalM(f_{1,2})$ and $\CalM(f_{1,5})$ which depends on $g$. Make the decomposition $g = g_0 + g_1$ with $g_0, g_1$ satisfying
\begin{align} \label{sys:g-0}
\begin{split}
 &  \del_t g_0 + \mu v_\omega \del_x g_0 
= - \frac{1}{\tau(\omega)} g_0, 
\\
& g_0 |_{x=0} = \eta_2 f(t, 1, \mu, \omega) + \zeta_2 g(t, 1, -\mu, \omega), \hspace{2.7cm} \mu > 0,   
\\
&  g_0|_{x=L} = \alpha_0 \xi(\omega) \int_0^1 \int_0^\infty \mu v_\omega g(t, L, \mu, \omega) \domega\dmu, \hspace{1.4cm} \mu < 0,  
\\
& g_0 |_{t=0} = 0, 
\end{split}
\end{align}
and
\begin{align} \label{sys:g-1}
\begin{split}
 &  \del_t g_1 + \mu v_\omega \del_x g_1 
= - \frac{1}{\tau(\omega)} g_1 + \frac{\xi(\omega)}{\tau(\omega)} \dvint{g}, 
\\
& g_1 |_{x=0} = 0, \hspace{3cm} \mu > 0,   
\\
&  g_1|_{x=L} = 0, \hspace{3cm} \mu < 0,  
\\
& g_1 |_{t=0} = 0. 
\end{split}
\end{align}
First we show that with assumption (A4), $g_0$ does not contribute to either $\CalM(f_{1, 2})$ or $\CalM(f_{1, 5})$. This is because its support (in time) is outside of the support of $\psi(\frac{t-t_1}{\Eps})$. 
More specifically, we start by noting that the boundary data of $g_0$ do not play a role. In fact, if we denote the data of $g_0$ at $x=1$ as $H_1(t, \mu, \omega) \One_{\mu > 0}$ and at $x=L$ as $H_L(t, \mu, \omega) \One_{\mu < 0}$, then by solving along the characteristics, we obtain that
\begin{align*}
   H_L(t, \mu, \omega) \One_{\mu < 0} = 0,
\qquad
  0 \leq t \leq \tfrac{L-1}{|\mu| v_\omega}.
\end{align*}
Hence, 
\begin{align*}
   g_0(t, 1, \mu, \omega) \One_{\mu<0} = 0, 
\qquad
  0 \leq t \leq \tfrac{2(L-1)}{|\mu| v_\omega}.
\end{align*}
This gives 
\begin{align*}
   g_0(t + \tfrac{1}{\mu v_\omega}, 1, \mu, \omega) \One_{\mu < 0} = 0, 
\qquad
  \tfrac{1}{|\mu| v_\omega} \leq t \leq \tfrac{2L-1}{|\mu| v_\omega}.
\end{align*}
By~\eqref{assump:L} in assumption (A4), we have
\begin{align*}
   \tfrac{2L-1}{|\mu| v_\omega}
\geq 
  \tfrac{2L-1}{v_0}
> \tfrac{2}{\mu_0 v_{\omega_0}} + 1
\geq 
   t_1 + \Eps,
\end{align*}
which shows 
\begin{align*}
   g_0(t + \tfrac{1}{\mu v_\omega}, 1, \mu, \omega) = 0, 
\qquad
   \forall t \in \text{Supp} \, \psi(\tfrac{t-t_1}{\Eps}). 
\end{align*}
Hence, the contribution of $g_0$ to $\CalM(f_{1,2}) + \CalM(f_{1,5})$ is zero.

Finally, we estimate the contribution of~$g_1$, which comes from the inhomogeneous term $\dvint{g}$.  The estimate will again follow from Proposition~\ref{prop:general}.
Directly solving the $g_1$-equation gives
\begin{align*}
   g_1(t, 1, \mu, \omega)
= \begin{cases}
   \int_0^t e^{-\frac{y}{\tau(\omega)}}
  \tfrac{\xi(\omega)}{\tau(\omega)}
  \dvint{g}(t - y, 1 - \mu v_\omega y) \dy, 
\qquad
  \mu < 0,
\quad
  0 < t < \tfrac{L-1}{|\mu| v_\omega},  \\[4pt]
  \int_0^{\frac{L-1}{|\mu| v_\omega}} e^{-\frac{y}{\tau(\omega)}}
  \tfrac{\xi(\omega)}{\tau(\omega)}
  \dvint{g}(t - y, 1 - \mu v_\omega y) \dy, 
\qquad
  \mu < 0,
\quad
  t > \tfrac{L-1}{|\mu| v_\omega}. 
\end{cases}
\end{align*}
Thus in either cases, we have
\begin{align*}
   g_1(t, 1, \mu, \omega)
\leq \int_0^t e^{-\frac{y}{\tau(\omega)}}
  \tfrac{\xi(\omega)}{\tau(\omega)}
  \dvint{g}(t - y, 1 - \mu v_\omega y) \dy.
\end{align*}
Combining $\CalM(f_{1,2})$ and $\CalM(f_{1,5})$ with $g$ replaced by $g_1$, we have
\begin{align*}
   \abs{\CalM(f_{1,2}+f_{1,5})}
&= \abs{\int_0^\infty \int_{-1}^0 \int_0^\infty 
    \psi_0\vpran{\tfrac{t-t_1}{\Eps}}
    \zeta_1(\omega) g_1 \vpran{t + \tfrac{1}{\mu v_\omega}, 1, \mu, \omega} e^{-\frac{1}{\tau(\omega)} \frac{1}{|\mu| v_\omega}} 
    \One_{t \geq \frac{1}{|\mu| v_\omega}} \dt\dmu\domega}
\\
& \leq
   C \int_0^\infty \int_{-1}^0 \int_0^\infty
   \int_0^{t + \frac{1}{\mu v_\omega}}
   \psi_0\vpran{\tfrac{t-t_1}{\Eps}} \frac{\xi(\omega)}{\tau(\omega)}
   \dvint{|g|}\vpran{t + \tfrac{1}{\mu v_\omega} - y, 1 - \mu v_\omega y}
   \One_{t \geq \frac{1}{|\mu| v_\omega}} \dy\dt\dmu\domega
\end{align*}
which has a similar structure as $\CalM(f_{1,3})$. Hence, 
\begin{align} \label{bound:M-f-1-part-2}
  \abs{\CalM(f_{1,2} + f_{1,5})} 
\leq 
  C_{p, \phi_0} \Eps^{1-\frac{3}{p'_0}} \to 0
\qquad
  \text{as}\quad \Eps \to 0\,.
\end{align}
Combining~\eqref{bound:M-f-1-part-1} with~\eqref{bound:M-f-1-part-2} we obtain the desired bound in~\eqref{bound:M-f-1-total}. 
\end{proof}

As the concluding remark, the proof of Theorem~\ref{thm:main} follows from Propositions~\ref{prop:f0} and~\ref{prop:f1}. 

\section{Conclusion}
We prove the unique reconstruction of the heat-reflection indices at the interface of two solids in an experiment setup presented in~\cite{PhysRevB.95.205423}. Our result thus justifies the validity of this experiment. For the rigour of the proof, we use highly singular data concentrated around $\mu_0$, $\omega_0$ and $t_0$. This triggers the response of the heat-reflection indices with frequency $\omega_0$ at a particular time and the measurement should be placed at $t_0+\frac{2}{\mu_0v_{\omega_0}}$ to record the reflected heat, which reveals the heat indices at $\omega_0$. In experiments, however, it is hard to prepare an input that concentrates in time -- the prescribed incoming data usually oscillate in time with a preset frequency (dual of time). This is not a contradiction to our theory, since by the linearity of the equation, the oscillatory incoming data can be viewed as a superposition of singular incoming data. We leave the investigation of the real experimental data to the future research.


\bibliographystyle{amsxport}
\bibliography{transbib}

\end{document}